\documentclass[a4paper,12pt]{article}

 \usepackage{latexsym,amssymb,amsmath}

 \newtheorem{prop}{Proposition}

\begin{document}

\title{Vekua-type systems related to two-sided monogenic functions}

\author{Dixan Pe\~na Pe\~na$^{\star,1}$ and Frank Sommen$^{\star\star,2}$}

\date{\normalsize{$^\star$Department of Mathematics, Aveiro University\\3810-193 Aveiro, Portugal\\
$^{\star\star}$Department of Mathematical Analysis, Ghent University\\9000 Ghent, Belgium}\\\vspace{0.4cm}
\small{$^1$e-mail: dixanpena@ua.pt; dixanpena@gmail.com\\
$^2$e-mail: fs@cage.ugent.be}} 

\maketitle

\begin{abstract}
\noindent Solutions to the Dirac equation are obtained by considering functions of axial type. This indeed gives rise to Vekua systems that can be solved in terms of special functions. In this paper we investigate axial symmetry for the solutions of the two-sided monogenic system and we give examples involving Bessel functions.\vspace{0.2cm}\vspace{0.1cm}
\end{abstract}

\section{Introduction}

Let $\mathbb{R}_{0,m}$ be the $2^m$-dimensional real Clifford algebra constructed over the orthonormal basis $(e_1,\ldots,e_m)$ of the Euclidean space $\mathbb R^m$ (see \cite{Cl}). The multiplication in $\mathbb{R}_{0,m}$ is governed by the following rules 
\begin{alignat*}{2} 
e_j^2=-1&,&\qquad &j=1,\dots,m,\\
e_je_k+e_ke_j=0&,&\qquad &1\le j\neq k\le m.
\end{alignat*}
A basis for the algebra $\mathbb{R}_{0,m}$ is given by the elements 
\[e_A=e_{j_1}\cdots e_{j_k},\]
where $A=\{j_1,\dots,j_k\}\subset\{1,\dots,m\}$ is such that $j_1<\dots<j_k$. For the empty set $\emptyset$, we put $e_{\emptyset}=1$, the latter being the identity element. The subspace of $k$-vectors is defined as
\[\mathbb R_{0,m}^{(k)}=\Bigg\{a\in\mathbb R_{0,m}:\;a=\sum_{\vert A\vert=k}a_Ae_A\Bigg\},\]
leading to the decomposition $\mathbb R_{0,m}=\bigoplus_{k=0}^m\mathbb R_{0,m}^{(k)}$.

One way to generalize the holomorphic functions of a complex variable is by considering the null solutions of the so-called generalized Cauchy-Riemann operator in $\mathbb R^{m+1}$, given by
\[\partial_x=\partial_{x_0}+\partial_{\underline x},\]
where 
\[\partial_{\underline x}=\sum_{j=1}^me_j\partial_{x_j}\] 
is the Dirac operator in $\mathbb R^m$. More precisely, an $\mathbb R_{0,m}$-valued function $f$ defined and continuously differentiable in an open set $\Omega$ of $\mathbb R^{m+1}$, is said to be left (resp. right) monogenic in $\Omega$ if $\partial_xf=0$ (resp. $f\partial_x=0$) in $\Omega$ (see e.g. \cite{BDS,DSS,GuSp}). 

In a similar fashion is defined monogenicity with respect to $\partial_{\underline x}$. Additionally, functions which are both left and right monogenic are called two-sided monogenic.

Suppose that $P_k(\underline x)$ is a homogeneous left monogenic polynomial of degree $k\in\mathbb N_0$ in $\mathbb R^m$. The axial (left) monogenic functions are monogenic functions of the form 
\begin{equation}\label{AxialLMF}
\Bigl(A(x_0,r)+\frac{\underline x}{r}\,B(x_0,r)\Bigr)P_k(\underline x),\quad r=\vert\underline x\vert=\sqrt{-\underline x^2},
\end{equation}
with $A$ and $B$ being $\mathbb R$-valued and continuously differentiable functions in $\mathbb R^2$ (see \cite{LB,S1,S2}). It is not difficult to show that $A$ and $B$ satisfy the Vekua-type system
\begin{equation*}\label{Veq}
\left\{\begin{array}{ll}\partial_{x_0}A-\partial_rB&=\displaystyle{\frac{2k+m-1}{r}}\,B\\\partial_{x_0}B+\partial_rA&=0.\end{array}\right.
\end{equation*}
It is worth pointing out that some examples of these functions are provided by the so-called Fueter's theorem (see e.g. \cite{F,KQS,DS1,DS2,Sce,S3}). 

Our aim in this paper is to construct the analog of functions (\ref{AxialLMF}) for the case of two-sided monogenicity. Explicit examples will also be given.

\section{Two-sided axial monogenic functions}

Let us recall the following Leibniz rules that will be used in our calculations
\begin{equation}\label{lr1}
\partial_{\underline x}(\underline xf)=-mf-2\sum_{j=1}^mx_j(\partial_{x_j}f)-\underline x(\partial_{\underline x}f),
\end{equation}
\begin{equation*}\label{lr2}
(f\underline x)\partial_{\underline x}=-mf-2\sum_{j=1}^mx_j(\partial_{x_j}f)-(f\partial_{\underline x})\underline x.
\end{equation*}
We denote by $P_{k,l}(\underline x)$ a homogeneous two-sided monogenic polynomial of degree $k$ in $\mathbb R^m$ with values in the subspace of $l$-vectors $\mathbb R_{0,m}^{(l)}$. To be precise:
\begin{alignat*}{2}
P_{k,l}(t\underline x)&=t^kP_{k,l}(\underline x),&\quad t&\in\mathbb R,\\
\partial_{\underline x}P_{k,l}(\underline x)&=0=P_{k,l}(\underline x)\partial_{\underline x},&\quad\underline x&\in\mathbb R^m,\\
P_{k,l}(\underline x)&\in\mathbb R_{0,m}^{(l)},&\quad\underline x&\in\mathbb R^m.
\end{alignat*}
In what follows, $J_\alpha$, $Y_\alpha$ stand for the Bessel functions of the first and second kind respectively (see e.g. \cite{Hoch}). They satisfy the recurrence relations:  
\begin{align*}
\frac{2\alpha}{t}Z_\alpha(t)&=Z_{\alpha-1}(t)+Z_{\alpha+1}(t),\\
2Z_\alpha^\prime(t)&=Z_{\alpha-1}(t)-Z_{\alpha+1}(t),
\end{align*}
where $Z_\alpha$ denotes $J_\alpha$ or $Y_\alpha$. On account of the above relations, we have
\[\left(t^{-\alpha}Z_\alpha(t)\right)^\prime=-t^{-\alpha}Z_{\alpha+1}(t).\]
Let us now assume that $A$, $B$, $C$, $D$ are $\mathbb R$-valued continuously differentiable functions in some open subset of $\mathbb R^2_+=\{(x_1,x_2)\in\mathbb R^2:\;x_2>0\}$, and consider
\begin{multline*}
F(x)=A(x_0,r)P_{k,l}(\underline x)+B(x_0,r)\underline xP_{k,l}(\underline x)\\+C(x_0,r)P_{k,l}(\underline x)\underline x+D(x_0,r)\underline x P_{k,l}(\underline x)\underline x.
\end{multline*}
It is easily seen that 
\[\partial_{\underline x}A=\sum_{j=1}^me_j\partial_{x_j}A=\sum_{j=1}^me_j(\partial_rA)(\partial_{x_j}r)=\frac{\partial_rA}{r}\,\underline x\]
and therefore
\[\partial_{\underline x}(AP_{k,l})=(\partial_{\underline x}A)P_{k,l}+A(\partial_{\underline x}P_{k,l})=\frac{\partial_rA}{r}\underline xP_{k,l}.\]
Using the Leibniz rule (\ref{lr1}) and Euler's theorem for homogeneous functions, we also obtain that
\begin{align*}
\partial_{\underline x}(B\underline xP_{k,l})&=(\partial_rB)\frac{\underline x^2}{r}P_{k,l}-B\big(mP_{k,l}+2\sum_{j=1}^mx_j(\partial_{x_j}P_{k,l})+\underline x(\partial_{\underline x}P_{k,l})\big)\\
&=-\big((2k+m)B+r\partial_rB\big)P_{k,l}.
\end{align*}
Since, by hypothesis $P_{k,l}$ is an $l$-vector valued function, then
\begin{equation*}
\sum_{j=1}^me_jP_{k,l}e_j=(-1)^l(2l-m)P_{k,l},
\end{equation*}
which leads to
\[\partial_{\underline x}(P_{k,l}\underline x)=(\partial_{\underline x}P_{k,l})\underline x+\sum_{j=1}^me_jP_{k,l}(\partial_{x_j}\underline x)=(-1)^l(2l-m)P_{k,l}.\]
This gives
\begin{align*}
\partial_{\underline x}(CP_{k,l}\underline x)&=(-1)^l(2l-m)CP_{k,l}+\frac{\partial_rC}{r}\underline xP_{k,l}\underline x,\\
\partial_{\underline x}(D\underline xP_{k,l}\underline x)&=(-1)^{l+1}(2l-m)D\underline xP_{k,l}\\
&\qquad\qquad\qquad\quad-\big((2k+m+2)D+r\partial_rD\big)P_{k,l}\underline x.
\end{align*}
In view of the above equalities, we see that $F$ is left monogenic if 
\begin{equation}\label{Veq2sidedI}
\left\{\begin{aligned}
\partial_{x_0}A-r\partial_rB&=(2k+m)B+(-1)^{l+1}(2l-m)C\\
\partial_{x_0}B+\frac{1}{r}\partial_rA&=(-1)^l(2l-m)D\\
\partial_{x_0}C-r\partial_rD&=(2k+m+2)D\\
\partial_{x_0}D+\frac{1}{r}\partial_rC&=0.
\end{aligned}\right.
\end{equation}
Similarly, we can also get that 
\[\hspace{-8.2cm}(AP_{k,l})\partial_{\underline x}=\frac{\partial_rA}{r}P_{k,l}\underline x,\]
\[\hspace{-4.2cm}(B\underline xP_{k,l})\partial_{\underline x}=(-1)^l(2l-m)BP_{k,l}+\frac{\partial_rB}{r}\underline xP_{k,l}\underline x,\]
\[\hspace{-5.3cm}(CP_{k,l}\underline x)\partial_{\underline x}=-\big((2k+m)C+r\partial_rC\big)P_{k,l},\]
\[(D\underline xP_{k,l}\underline x)\partial_{\underline x}=-\big((2k+m+2)D+r\partial_rD\big)\underline xP_{k,l}+(-1)^{l+1}(2l-m)DP_{k,l}\underline x,\]
leading to the following system for the case of right monogenicity of $F$
\begin{equation}\label{Veq2sidedII}
\left\{\begin{aligned}
\partial_{x_0}A-r\partial_rC&=(2k+m)C+(-1)^{l+1}(2l-m)B\\
\partial_{x_0}C+\frac{1}{r}\partial_rA&=(-1)^l(2l-m)D\\
\partial_{x_0}B-r\partial_rD&=(2k+m+2)D\\
\partial_{x_0}D+\frac{1}{r}\partial_rB&=0.
\end{aligned}\right.
\end{equation}
Since we are interested if $F$ being two-sided monogenic, we should consider simultaneous solutions of the systems (\ref{Veq2sidedI})-(\ref{Veq2sidedII}). A quick look at these systems reveals that 
\[\partial_{x_0}(B-C)=\partial_r(B-C)=0,\]
which clearly implies that $B-C$ is constant. Subtracting the first equations from (\ref{Veq2sidedI})-(\ref{Veq2sidedII}) we may conclude that $B=C$.  

\begin{prop}
Suppose that $A_1$, $A_2$, $A_3$ are $\mathbb R$-valued continuously differentiable functions in the open set $\Xi\subset\mathbb R^2_+$. If these functions satisfy in $\Xi$ the overdetermined system 
\begin{equation}\label{Veq2sidedI-II}
\left\{\begin{aligned}
\partial_{x_0}A_1-r\partial_rA_2&=\big((2k+m)+(-1)^{l+1}(2l-m)\big)A_2\\
\partial_{x_0}A_2+\frac{1}{r}\partial_rA_1&=(-1)^l(2l-m)A_3\\
\partial_{x_0}A_2-r\partial_rA_3&=(2k+m+2)A_3\\
\partial_{x_0}A_3+\frac{1}{r}\partial_rA_2&=0,
\end{aligned}\right.
\end{equation}
then 
\[A_1(x_0,r)P_{k,l}(\underline x)+A_2(x_0,r)\underline xP_{k,l}(\underline x)+A_2(x_0,r)P_{k,l}(\underline x)\underline x+A_3(x_0,r)\underline x P_{k,l}(\underline x)\underline x\]
is two-sided monogenic in $\Xi^*=\{x\in\mathbb R^{m+1}:\;(x_0,r)\in\Xi\}$. 
\end{prop}

We now try to find particular solutions of the system (\ref{Veq2sidedI-II}), which we assume to be of the form 
\[A_j(x_0,r)=\exp(x_0)a_j(r),\quad j=1,2,3,\] 
with $a_1$, $a_2$, $a_3$ being $\mathbb R$-valued continuously differentiable functions.

From (\ref{Veq2sidedI-II}) it follows that
\begin{equation}\label{Veq2sidedIII}
\left\{\begin{aligned}
a_1-ra_2^\prime&=\big((2k+m)+(-1)^{l+1}(2l-m)\big)a_2\\
a_2+\frac{a_1^\prime}{r}&=(-1)^l(2l-m)a_3\\
a_2-ra_3^\prime&=(2k+m+2)a_3\\
a_3+\frac{a_2^\prime}{r}&=0.
\end{aligned}\right.
\end{equation}
Eliminating $a_3$ from the last two equations of (\ref{Veq2sidedIII}), yields
\[ra_2^{\prime\prime}+(2k+m+1)a_2^\prime+ra_2=0.\]
The general solution of this homogeneous ordinary differential equation is expressed in terms of the Bessel functions: 
\begin{equation*}
a_2(r)=r^{-k-\frac{m}{2}}\left(C_1J_{k+\frac{m}{2}}(r)+C_2Y_{k+\frac{m}{2}}(r)\right),
\end{equation*}
where $C_1$, $C_2$ are arbitrary real constants. Hence
\begin{equation*}
a_3(r)=r^{-k-\frac{m}{2}-1}\left(C_1J_{k+\frac{m}{2}+1}(r)+C_2Y_{k+\frac{m}{2}+1}(r)\right).
\end{equation*}
We thus obtain, from the first equation of (\ref{Veq2sidedIII}), that
\[a_1(r)=\big((2k+m)+(-1)^{l+1}(2l-m)\big)a_2(r)-r^2a_3(r).\]
Finally, it is not difficult to check that $a_1$, $a_2$, $a_3$ fulfill the second equation of (\ref{Veq2sidedIII}).

\begin{prop}
Let $a_1$, $a_2$, $a_3$ be as above. Then the function 
\[\exp(x_0)\Big(a_1(r)P_{k,l}(\underline x)+a_2(r)\underline xP_{k,l}(\underline x)\\+a_2(r)P_{k,l}(\underline x)\underline x+a_3(r)\underline xP_{k,l}(\underline x)\underline x\Big)\]  
is two-sided monogenic in $\mathbb R^{m+1}\setminus\{\underline x=0\}$. 
\end{prop}

We would like to remark that solutions of the system (\ref{Veq2sidedI-II}) can also be obtained by using the power series method. Indeed, writing
\[A_j(x_0,r)=\sum_{n=0}^\infty\frac{x_0^n}{n!}A_{j,n}(r),\quad j=1,2,3,\]
and substituting into (\ref{Veq2sidedI-II}) yields the recurrence relations
\begin{align*}
A_{2,n+1}&=-A_{2,n-1}^{\prime\prime}-\frac{2k+m+1}{r}A_{2,n-1}^\prime,\\
A_{3,n+1}&=-\frac{1}{r}A_{2,n}^\prime,\\
A_{1,n+1}&=rA_{2,n}^\prime+\big((2k+m)+(-1)^{l+1}(2l-m)\big)A_{2,n},
\end{align*}
with initial conditions
\begin{align*}
A_{j,0}(r)&=A_j(0,r),\quad j=1,2,3,\\
A_{2,1}(r)&=\partial_{x_0}A_2(0,r).
\end{align*}
Clearly,
\begin{align*}
A_{2,2n}(r)&=\sum_{j=1}^{2n}c_{n,j}\,\frac{A_{2,0}^{(2n-j+1)}(r)}{r^{j-1}},\\
A_{2,2n+1}(r)&=\sum_{j=1}^{2n}c_{n,j}\,\frac{A_{2,1}^{(2n-j+1)}(r)}{r^{j-1}},
\end{align*}
with $n\in\mathbb N$ and $c_{n,j}\in\mathbb Z$.

\subsection*{Acknowledgment}

The first author was supported by a Post-Doctoral Grant of \emph{Funda\c{c}\~ao para a Ci\^encia e a Tecnologia}, Portugal (grant number: SFRH/BPD/45260/2008).

\end{document}